\documentclass[11pt]{amsart}
\usepackage{latexsym}
\usepackage{amssymb,amscd,amsthm,amsfonts,verbatim,amsopn,color}
\usepackage[dvips]{epsfig}
\def\R{{\mathbb R}}

\def\FF{{\mathcal F}}
\def\GG{{\mathcal G}}
\def\HH{{\mathcal H}}
\def\LL{{\mathcal L}}

\def\NN{{\mathcal N}}

\def\Ss{{\mathcal S}}

\newcommand{\lra}{\longrightarrow}

\newtheorem{thm}{Theorem}[section]
\newtheorem{df} [thm]{Definition}
\newtheorem{lm}  [thm]{Lemma}
\newtheorem{crl} [thm]{Corollary}
\newtheorem{prop}[thm]{Proposition}

\newtheorem{rem}[thm]{Remark}

\newtheorem{expl}[thm]{Example}
\numberwithin{equation}{section}

\newenvironment{pf}{\noindent {\bf Proof.}}{\hfill $\Box$\vspace{0.3cm}}

\newenvironment{remrm}{\begin{rem} \rm}{\end{rem}}
\newenvironment{dfrm}{\begin{df} \rm}{\end{df}}

\begin{document}

\title[On the topology of nested set complexes]
{On the topology of nested set complexes}

\author{Eva Maria Feichtner \mbox{ }\& \mbox{ }Irene M\"uller}

\address{
Department of Mathematics, ETH Zurich, 8092 Zurich, Switzerland}
\email{feichtne@math.ethz.ch} 
\email{irene@math.ethz.ch} 

\date{ 
November 2003. 
}

\begin{abstract}
  Nested set complexes appear as the combinatorial core of
  De~Concini-Procesi arrangement models. We show that nested set
  complexes are homotopy equivalent to the order complexes of the
  underlying meet-semilattices without their minimal elements.  For
  atomic semilattices, we consider the realization of nested set
  complexes by simplicial fans proposed in~\cite{FY}, and we
  strengthen our previous result showing that in this case nested set
  complexes in fact are homeomorphic to the mentioned order complexes.
\end{abstract}

\maketitle

\section{Introduction}
In the same way as intersection lattices 
capture the combinatorial essence  of arrangements of hyperplanes,
building sets and nested set complexes encode the combinatorics
of De~Concini-Procesi arrangement models: They prescribe the model 
construction by sequences of blowups, they describe the incidence 
combinatorics of the divisor stratification, and they naturally appear 
in presentations of cohomology algebras for arrangement models 
in terms of generators and relations (cf.\ \cite{DP}).

Nested set complexes have been defined in various generalities. 
The notion of nested sets goes back to the model construction
for configuration spaces of algebraic varieties by Fulton~\& 
MacPherson~\cite{FM}; the underlying poset in this special case 
is the lattice of set partitions. 
De~Concini and Procesi~\cite{DP} defined building 
sets and nested set complexes for intersection lattices of subspace 
arrangements in real or complex linear space; in this setting they 
have the broad geometric significance outlined above.

In joint work of the first author with D.\ Kozlov~\cite{FK},  
purely  order-theoretic definitions of building sets and nested set 
complexes for arbitrary meet-semilattices were given. 
Together with the notion of a combinatorial blowup in a meet-semilattice 
a complete combinatorial counterpart to the resolution process of 
DeConcini and Procesi was established. Having these purely combinatorial 
notions at hand, Yuzvinsky and the first author \cite{FY} studied abstract 
algebras that generalize arrangement model cohomology and solely depend on 
nested set complexes. In this context, nested set complexes attain yet
another geometric meaning as the defining data for certain toric 
varieties.

In this article, we study nested set complexes from the viewpoint of 
topological combinatorics. Relying on
techniques from the homotopy theory of partially ordered sets due 
to Quillen~\cite{Q}, we show 
that, for any building set $\GG$ in a meet-semilattice $\LL$, 
the nested set complex $\NN(\LL,\GG)$ is homotopy equivalent 
to the  order complex of the underlying meet-semilattice without its
minimal element, $\LL{\setminus}\{\hat 0\}$.

For atomic meet-semilattices we can strengthen this result. We consider
the realization of nested set complexes $\NN(\LL,\GG)$ by simplicial fans 
$\Sigma(\LL,\GG)$ proposed 
in~\cite{FY}, and we show that, for building sets $\HH\,{\subseteq}\,\GG$ 
in~$\LL$, the simplicial fan $\Sigma(\LL,\GG)$ is obtained from 
$\Sigma(\LL,\HH)$ 
by a sequence of stellar subdivisions. This in particular implies that,
for a given atomic meet-semilattice~$\LL$, 
the nested set complex for any building set is homeomorphic 
to the order complex of~$\LL{\setminus}\{\hat 0\}$. 

After a brief review of the definitions for building sets, nested set 
complexes, and combinatorial blowups in Section~\ref{sec_prel}, we present 
our result on the homotopy type of nested set complexes 
in Section~\ref{sec_htpy}.  
The strengthening in the case of atomic meet-semilattices is given in 
Section~\ref{sec_fans}.

{\em Acknowledgments:} We would like to thank Sergey Yuzvinsky 
for helpful discussions. 


\section{Preliminaries on building sets and nested sets}
\label{sec_prel}

For the sake of completeness we here review the definitions of building 
sets and nested sets for finite meet-semilattices as proposed in~\cite{FK}.

All posets occurring in this article are finite. We mostly assume that 
the posets are meet-semilattices (semilattices, for short), 
i.e., greatest lower bounds exist for any subset of elements in the poset. 
Any finite meet-semilattice $\LL$ has a minimal element, 
which we denote by~$\hat 0$. We frequently use the notation $\LL_{> \hat 0}$
to denote $\LL$ without its minimal element. 
For any subset $\Ss$ in $\LL$ we denote the 
set of maximal elements in $\Ss$ by $\mathrm{max}\, \Ss$. For any $X\,{\in}\,\LL$,
we set $\Ss_{\leq X}\,{=}\,\{Y\,{\in}\,\Ss\,|\,Y\,{\leq}\,X\}$, and we use the 
standard notation for intervals in~$\LL$, 
$[X,Y]\,{:=}\,\{Z\in \LL\,|\, X\,{\leq}\,Z\,{\leq}\,Y\}$. 
The standard simplicial complex built from a poset~$\LL$ is the 
{\em order complex\/} of $\LL$,
which we denote by $\Delta (\LL)$; it is the abstract simplicial complex on 
the elements of~$\LL$ with simplices corresponding to linearly ordered subsets 
in $\LL$. As a general reference on posets we refer to~\cite[Chapter~3]{St}.

\begin{dfrm} \label{df_building}
Let $\LL$ be a finite meet-semilattice. A subset $\GG$ 
in~$\LL_{>\hat 0}$  is called a {\em building set\/} if for any 
$X\,{\in}\,\LL_{>\hat 0}$  and
{\rm max}$\, \GG_{\leq X}=\{G_1,\ldots,G_k\}$ there is an isomorphism
of posets
\begin{equation}\label{eq_buildg}
\varphi_x:\,\,\, \prod_{j=1}^k\,\,\, [\hat 0,G_j] \,\, 
                               \stackrel{\cong}{\lra}
                                \,\, [\hat 0,X]
\end{equation}
with $\varphi_X(\hat 0, \ldots, G_j, \ldots, \hat 0)\, = \, G_j$ 
for $j=1,\ldots, k$. 
We call $F_{\GG}(X)\,{:=}\, 
\mathrm{max}\,\GG_{\leq X}$  the {\em set of factors\/} of $X$ in $\GG$.
\end{dfrm}

As a simple example we can take the full semilattice $\LL_{>\hat 0}$ 
as a building set. Besides this maximal building set, there is a minimal 
building set consisting of all elements $X$ in $\LL_{>\hat 0}$ which 
do not allow for a 
product decomposition of the lower interval~$[\hat 0,X]$, 
the so-called {\em irreducible elements\/} in~$\LL$.

Any choice of a building set $\GG$ in $\LL$ gives rise to a family of
so-called {\em nested sets\/}. These are, roughly speaking, subsets of
$\GG$ whose antichains are sets of factors with respect to the
building set $\GG$.  Nested sets form an abstract simplicial complex
on the vertex set $\GG$ -- the {\em nested set complex\/}, which is
the main character of this article.

\begin{dfrm} \label{df_nested}
Let $\LL$ be a finite meet-semilattice and $\GG$ a building set in $\LL$.
A subset $\Ss$ in $\GG$ is called~{\em nested\/} (or $\GG$-{\em nested\/}
if specification is needed) 
if, for any set of incomparable elements 
$X_1,\dots,X_t$ in $\Ss$ of cardinality at least two, 
the join $X_1\vee\dots\vee X_t$ exists and does not belong to $\GG$.
The $\GG$-nested sets form an abstract simplicial complex $\NN(\LL,\GG)$,
the {\em nested set complex\/} with respect to $\LL$ and $\GG$. 
\end{dfrm}

For the maximal building set $\LL_{>\hat 0}$ in $\LL$, the nested sets are the 
chains in $\LL_{>\hat 0}$; in particular, the nested set complex 
$\NN(\LL,\LL_{>\hat 0})$ coincides with the order complex 
$\Delta(\LL_{>\hat 0})$.

\medskip
We also remind here a construction on semilattices that was proposed 
in~\cite{FK}, the combinatorial blowup of a semilattice $\LL$ in an 
element  $X$ of $\LL$.

\begin{dfrm} \label{df_blowup}
For a semilattice $\LL$ and an element $X$ in $\LL_{>\hat 0}$ we define 
a poset $(\mathrm{Bl}_X\LL,\prec)$ on the set of elements
\[
  \mathrm{Bl}_X\LL\, \, = \, \,  
   \{\,Y\,|\, Y\in \LL, Y\not \geq X\} \,\, \cup \, \,  
   \{\,\hat Y\,|\, Y\in \LL, Y\not \geq X, \,\,\mathrm{and}\,\,
                   Y\vee X \,\, \mathrm{exists \, \, in}\,\,\LL \,\} \, .
\]
The order relation $<$ in $\LL$ determines the order relation 
$\prec$ within the two parts of~$\mathrm{Bl}_X\LL$ described above,
\[
\begin{array}{lll}
    Y \prec Z\, , & \mathrm{for} &  Y<Z \,\,\,\mathrm{in}\,\, \LL\, , \\
    \hat Y \prec \hat Z\, , & \mathrm{for} &  Y<Z \,\,\,\mathrm{in}\,\, \LL\, , \\
    \end{array}
\]
and additional order relations between elements of these two parts are 
defined by 
\[
\begin{array}{lll} 
  Y \prec \hat Z\, , & \mathrm{for} &  Y\leq Z \,\,\,\mathrm{in}\,\, \LL\, , 
\end{array}
\]
where in all three cases it is assumed that $Y,Z \not \geq X$ in $\LL$.
We call $\mathrm{Bl}_X\LL$ the {\em combinatorial blowup\/} of 
$\LL$ in $X$. 
\end{dfrm}

Let us remark here that $\mathrm{Bl}_X\LL$ is again a meet-semilattice. 
The combinatorial blowup of a semilattice was used in~\cite{FK} to analyze 
the incidence change of strata in the construction process for 
De~Concini-Procesi arrangement models. In the present paper we will
need combinatorial blowups to describe the incidence change in polyhedral
fans under stellar subdivision following an observation 
in~\cite[Prop.\ 4.9]{FK}:

\begin{prop}
Let $\Sigma$ be a polyhedral fan with face poset $\FF(\Sigma)$. For a 
cone $\sigma$ in $\Sigma$, the face poset of the fan obtained by stellar 
subdivision  of $\Sigma$ in $\sigma$, $\FF(\mathrm{st}(\Sigma,\sigma))$, 
can be described as the combinatorial blowup of $\FF(\Sigma)$ in $\sigma$:
\[
     \FF(\mathrm{st}(\Sigma,\sigma)) \, \, = \, \, 
     \mathrm{Bl}_{\sigma}(\FF(\Sigma)) \, .
\]
\end{prop}


\section{The homotopy type of nested set complexes} \label{sec_htpy}

In this section, we will show that for a given meet-semilattice $\LL$
and a building set $\GG$ in~$\LL$ the nested set complex $\NN(\LL,\GG)$ 
is homotopy 
equivalent to the order complex of~$\LL_{>\hat 0}$. We will use the 
following two lemmata on the homotopy type of partially ordered sets 
going back to Quillen~\cite{Q}.

\begin{lm} \label{lm_quillen}
{\rm (Quillen's fiber lemma)}
Let $f{:}\,P\,{\rightarrow}\,Q$ be a map of posets such that the order complex 
of $f^{-1}(Q_{\leq X})$ is contractible for all $X\,{\in}\,Q$, then
$f$ induces a homotopy equivalence between the order complexes of 
$P$ and~$Q$.  
\end{lm}

\begin{lm} \label{lm_joincontr}
Let $P$ be a poset, and assume that there is an element $X_0$ in $P$ 
such that the join $X_0\,{\vee}\,X$ exists for all $X\,{\in}\,P$. 
Then the order complex of $P$ is contractible. A poset with the property 
described above is called {\em join-contractible\/} via $X_0$.
\end{lm}

\begin{prop} \label{prop_htpy}
Let $\LL$ be a finite meet-semilattice, and $\GG$ a building set in $\LL$. 
Then the nested set complex $\NN(\LL,\GG)$ is homotopy equivalent to the 
order complex of $\LL_{>\hat 0}$,
\[
\NN(\LL,\GG)\, \, \simeq \, \, \Delta(\LL_{>\hat 0})\, .
\] 
\end{prop}

\begin{pf}
We denote by $\FF(\NN)$ the poset of non-empty faces of the nested set 
complex $\NN(\LL,\GG)$. Consider the following map of posets: 
\begin{eqnarray*}
    \phi\,: \, \, \FF(\NN) & \longrightarrow & \LL_{>\hat 0} \\
                     \Ss   & \longmapsto     & \bigvee\, \Ss\, = \, 
                                \bigvee_{X\in \Ss}\, X\, . 
\end{eqnarray*}
We claim that the order complex of 
$\FF(\NN)_{\leq X}\,{:=}\,\phi^{-1}((\LL_{>\hat 0})_{\leq X})$
is contractible for any $X\,{\in}\,\LL_{>\hat 0}$. An application of the 
Quillen fiber lemma~\ref{lm_quillen} will then prove the statement of the 
proposition. 

\noindent
{\em Case 1:\/} $X\,{\in}\,\GG$. We show that $\FF(\NN)_{\leq X}$ is 
join-contractible via $X$ and, with an application of 
Lemma~\ref{lm_joincontr}, thus prove our claim.
Let $\Ss$ be an element in $\FF(\NN)_{\leq X}$, i.e., a nested set with
$\bigvee \Ss\,{\leq}\,X$. We have to show that $\Ss\,{\cup}\,\{X\}$ is
nested with $\bigvee \Ss\,{\cup}\,\{X\}\, {\leq}\, X$, hence 
$\Ss\,{\cup}\,\{X\}\,{\in}\, \FF(\NN)_{\leq X}$.
Either $\bigvee \Ss\,{=}\,X$, in which case 
$X\,{\in}\,\Ss$, and our claim is obvious;
or $\bigvee \Ss\,{<}\,X$, in which case we can add $X$ to $\Ss$,
obtaining a nested set, with  $\bigvee \Ss\,{\cup}\,\{X\}\, {=}\, X$, 
hence $\Ss\,{\cup}\,\{X\}\,{\in}\, \FF(\NN)_{\leq X}$.

\noindent
{\em Case 2:\/} $X\,{\not\in}\,\GG$. We show that $\FF(\NN)_{\leq X}$ is 
join-contractible via the set of factors of~$X$, $F_{\GG}(X)$. 
Again, let $\Ss$ be a nested set with $\bigvee \Ss\,{\leq}\,X$;
we have to show that $\Ss\,{\cup}\,F_{\GG}(X)$ is nested 
with join less or equal~$X$, hence 
$\Ss\,{\cup}\,F_{\GG}(X)\,{\in}\,\FF(\NN)_{\leq X}$.

If  $\bigvee \Ss\,{=}\,X$, then $X\,{=}\,\bigvee {\rm max}\, \Ss$ and 
$F_{\GG}(X)\,{=}\, {\rm max}\,\Ss\,{\subseteq}\,\Ss$ by~\cite[Prop.\,2.8(2)]{FK},
which makes our claim obvious. 

For $\bigvee \Ss\,{<}\,X$, 
assume that $A\,{\subseteq}\,\Ss\,{\cup}\,F_{\GG}(X)$ is an antichain with at least
two elements, and $\bigvee A\,{\in}\,\GG$. Since the $\GG$-factors of $X$,
$F_{\GG}(X)\,{=}
\,\{G_1,\ldots,G_t\}$,
give a partition of $\GG_{\leq X}$ into subsets $\GG_{\leq G_i}$, 
$i\,{=}\,1,\ldots,t$~\cite[Prop.\,2.5(1)]{FK}, 
we find that $\bigvee A\,{\leq}\,G$ for some
$G\,{\in}\,F_{\GG}(X)$. 
If $A$ contains any elements of $F_{\GG}(X)$, then it must contain $G$,
which contradicts $A$ being an antichain with more than one element.
We conclude that $A$ does not contain any factors of $X$. In particular, 
it is a subset of the nested set $\Ss$, thus should have a join outside $\GG$,
and we again reach a contradiction. We conclude that $\Ss\,{\cup}\,F_{\GG}(X)$ 
is nested with join $X$, hence belongs to $\FF(\NN)_{\leq X}$. 
\end{pf}

\begin{remrm} \label{rm_crosscut}
The homotopy equivalence in \ref{prop_htpy} can be viewed 
as a generalization of the 
classical crosscut theorem~\cite[Thm.\,10.8]{Bj} applied to a particular class
of posets and crosscuts: Let $P$ be a simplicial poset, i.e., $P$ contains 
a minimal element $\hat 0$, and each interval $[\hat 0,X]$, $X\,{\in}\,P$,
is isomorphic to a boolean lattice. Observe that $P$ is a meet-semilattice,
and the set of atoms $\mathfrak A$ is a building set in $P$. 
The crosscut complex of $P$ with respect to $\mathfrak A$,
\[
    \Gamma(P,\mathfrak A)\, \, = \, \, 
      \{A\,{\subseteq}\,\mathfrak A\,|\, A \mbox{ is bounded in }\, P \}\, ,
\]
coincides with the nested set complex $\NN(P,\mathfrak A)$.
\end{remrm}


\section{Simplicial fans realizing nested set complexes} \label{sec_fans}

We recall the definition of the simplicial fan $\Sigma(\LL,\GG)$ for a given
atomic meet-semilattice~$\LL$ and a building set $\GG$ in $\LL$. 
For details see~\cite[Section 5]{FY}.  

Given a finite meet-semilattice $\LL$ with set of atoms
$\mathfrak{A}(\LL)\,{=}\,\{A_1,\ldots, A_n\}$, we will frequently use the
following notation: For $X\,{\in}\,\LL$, define 
$\lfloor X \rfloor\,{:=}\,\{A\,{\in}\,\mathfrak{A}(\LL)\, 
                                      |\, X\,{\geq}\,A\}$, 
the set of atoms below a specific element~$X$ in $\LL$.
We define characteristic
vectors $v_X$ in $\R^n$ for lattice elements 
$X\,{\in}\,\LL$ by
\[
(v_X)_i \, \, := \, \, \left\{ 
\begin{array}{ll}
1 & \mbox{ if }\, A_i \in \lfloor X \rfloor, \\
0 & \mbox{ otherwise}, \qquad \,\,\, \mbox{ for }\, i=1,\ldots,n.
\end{array}
\right. 
\]
These characteristic vectors will appear as spanning vectors of simplicial 
cones in $\R^n$. For a subset $\mathcal{S}\,{\subseteq}\,\LL$, we agree to 
denote by $V(\mathcal{S})$ the cone spanned by the vectors $v_X$ for 
$X\in \mathcal{S}$. 

\begin{dfrm}
Let~$\LL$ be a finite atomic meet-semilattice and~$\GG$ a building set
in~$\LL$. We define a rational, polyhedral fan $\Sigma(\LL,\GG)$ in $\R^n$ 
as the collection of cones $V(\Ss)$  for all nested sets  
$\Ss$ in $\LL$,
\begin{equation} \label{def_Sigma}
 \Sigma(\LL,\GG) \, \, :=\, \, 
             \{\, V(\mathcal{S})\,|\, \mathcal{S}\in \NN(\LL,\GG)\,\} \, . 
\end{equation}  
\end{dfrm}

By definition, rays in $\Sigma(\LL,\GG)$ are in $1$-$1$ correspondence 
with elements in~$\GG$. In fact, 
the face poset of $\Sigma(\LL,\GG)$ coincides with the face poset of 
$\NN(\LL,\GG)$. 

If there is no risk of confusion we will denote the fan in (\ref{def_Sigma})
by $\Sigma(\GG)$.

\begin{thm} \label{thm_stsubd}
Let $\LL$ be a finite atomic meet-semilattice, and $\GG$, $\HH$ building
sets in $\LL$ with $\GG\,{\supseteq}\,\HH$. Then, the fan $\Sigma(\GG)$
is obtained from $\Sigma(\HH)$ by a sequence of stellar subdivisions.
In particular, the supports of the fans  $\Sigma(\GG)$ and $\Sigma(\HH)$
coincide.
\end{thm}

\begin{pf}
For building sets $\GG\,{\supseteq}\,\HH$ in $\LL$ and $G$ minimal in
$\GG\,{\setminus}\,\HH$, set $\overline \GG\,{:=}\,\GG\,{\setminus}\,\{G\}$.
Obviously, ${\rm max}\,\overline \GG_{\leq G}\,{=}\,F_{\HH}(G)$, and for 
any $X\,{\in}\,\LL$ we find that
\[
     {\rm max} \, \overline \GG_{\leq X} \, = \, 
\left\{ 
\begin{array}{ll}
F_{\GG}(X) & \mbox{if }\, G\not\in F_{\GG}(X)\, , \\
(F_{\GG}(X)\setminus \{G\}) \cup F_{\HH}(G)
           & \mbox{if }\, G \in F_{\GG}(X)\, . \\
\end{array}
\right.
\]
Isomorphisms of posets required for the building set property of 
$\overline \GG$ expand accordingly in the second case, and we find that 
$\overline \GG$ is again a building set for~$\LL$. 

We thus conclude that, for any two building sets $\GG$, $\HH$
with $\GG\,{\supseteq}\,\HH$, there is a sequence of building sets
\[
   \GG \,= \, \GG_1 
       \, \supseteq \, \GG_2 
       \, \supseteq \, \ldots 
       \, \supseteq \, \GG_t
       \,=\,\HH\,  ,   
\]
such that $\GG_i$ and  $\GG_{i+1}$ differ by exactly one element $G_i$,
and $G_i$ is minimal in $\GG_i\,{\setminus}\,\HH$ for $i\,{=}\,1,\ldots,t{-}1$.

We can thus assume that $\HH\,{=}\,\GG\setminus\{G\}$, and it suffices 
to show that $\Sigma(\GG)$ is obtained from $\Sigma(\HH)$ by a sequence 
of stellar subdivisions.

In fact, we claim that $\Sigma(\GG)$ is obtained by a single stellar
subdivision of $\Sigma(\HH)$ in $V(F_{\HH}(G))$, introducing a new ray that is 
generated by the characteristic vector $v_G$, i.e.,
\begin{equation} \label{eq_fans}
\Sigma(\GG) \, \, =\, \, {\rm st}\,(\Sigma(\HH),V(F_{\HH}(G)),v_G )\, .
\end{equation}

Observe that the two fans in (\ref{eq_fans}) share the same set of
generating vectors for rays, so all we have to show is that they have the
same combinatorial structure, i.e., their face posets coincide.

The face poset of the subdivided fan can be described as the combinatorial 
blowup of the face poset $\FF(\Sigma(\HH))=\FF(\NN(\HH))$ in 
the $\HH$-nested set $F_{\HH}(G)$ (cf.~\cite[Sect.\ 4.2]{FK}), 
hence we are left to show that 
\begin{equation} \label{eq_posets}
\FF(\NN(\GG)) \, \, =\, \, {\rm Bl}_{F_{\HH}(G)}(\FF(\NN(\HH)))\, .
\end{equation}
Let us abbreviate notation and denote the poset on the right hand side
by $\mathrm{Bl}\,\FF$.

We first show the left-to-right inclusion in (\ref{eq_posets}). \newline
Let $\Ss$ be a $\GG$-nested set in $\LL$. We need to show that $\Ss$ is 
an element in $\mathrm{Bl}\,\FF$. For the matter of this proof, we agree 
to freely switch between sets of atoms and their joins in the respective 
semilattices.

For $G\,{\not\in}\,\Ss$, we note that $\Ss$ is $\HH$-nested. Moreover, 
$\Ss$ does not
contain $F_{\HH}(G)$, since the latter is certainly not $\GG$-nested.
We conclude that $\Ss$ is an element in~$\mathrm{Bl}\,\FF$.  

For $G\,{\in}\,\Ss$, we need to show that $(S\,{\setminus}\,\{G\})\,{\cup}\,
F_{\HH}(G)$ is $\HH$-nested. Let $A$ be an antichain with at least two elements 
that is contained in
$(S\,{\setminus}\,\{G\})\,{\cup}\, F_{\HH}(G)$; we need to see that 
$\bigvee A{\,\not\in}\,\HH$. If $A\,\subseteq\,F_{\HH}(G)$ then clearly 
$\bigvee A$ either equals $G$ or lives between $G$ and its $\HH$-factors, 
hence in any case is not contained in $\HH$. If $A$ does not contain
any $\HH$-factor of $G$, then $A\subset \Ss$ is $\GG$-nested, in particular 
$\bigvee A\,{\not\in}\,\HH$.

We can thus assume that the antichain $A$ is of the form
$A=\{S_1,\ldots, S_t,F_1,\ldots,F_k\}$, where 
$S_i\,{\in}\,\Ss\,{\setminus}\,(\,\{G\}\,{\cup}\, F_{\HH}(G)\,)$
for $i\,{=}\,1,\ldots,t$, and $F_j\,{\in}\, F_{\HH}(G)$  for $j\,{=}\,1,\ldots,k$,
and both types of elements occur in $A$.

Let us assume that $\bigvee A\,{\in}\,\HH$.  We have 
\begin{equation} \label{eq_A}
  \bigvee A \, \, \leq \, \, \bigvee_{i=1}^t S_i \, \,\, \vee \, \, \, G 
            \, \,   =  \, \, \bigvee_{
                 {i\in \{1,\ldots t\}, \, \,S_i \, \mathrm{ in-}} \atop
                                      {\mathrm{comparable \,\, with }\, G}}
                              S_i\,\,\, \vee \, \, \, G \,,
\end{equation}
where the last equality holds since any $S_j$ comparable with $G$ has to be
smaller than $G$, otherwise $S_j\,{\geq}\,G\,{>}\,F_1$ gives a contradiction 
to $A$ being an antichain. 

If there are no $S_i$, $i\,\in\,\{1,\ldots t\}$, that are incomparable with $G$,
the right hand side of (\ref{eq_A}) equals $G$. Assuming that 
$\bigvee A\,{\in}\,\HH$ we find that $\bigvee A\,{\leq}\, F$ for some 
$F\,{\in}\, F_{\HH}(G)$ since the $\HH$-factors of $G$ partition the elements 
of $\HH$ below $G$~\cite[Prop.\ 2.5.(1)]{FK}. 
We assumed that $A$ contains some of the $\HH$-factors of 
$G$, and thus conclude that it must contain $F$. This however contradicts 
to $A$ being an antichain with at least two elements.

We are left with the case of the join on the right hand side of (\ref{eq_A}) 
being taken over more than one element. Since 
$\Ss_0\,{=}\,\{S_i\,{\in}\,A\,|\, S_i$ incomparable with $G\}\,{\cup}\,\{G\} \, 
{\subseteq}\, \Ss$ is a $\GG$-nested antichain,  
we conclude that $\bigvee \Ss_0$ is not contained in $\GG$ and
$\Ss_0$ is its set of factors. Since these factors partition $\GG$-elements
below  $\bigvee \Ss_0$ we find that either
$\bigvee A \,{\leq}\, S_i$, for some $S_i\,{\in}\,\Ss_0$,
which is a contradiction to $A$ being an antichain, or 
$\bigvee A \,{\leq}\, G$, which again places $\bigvee A$ below one of 
the $\HH$-factors $F$ of $G$, and, as argued above, leads to a contradiction. 
We conclude that $(S\,{\setminus}\,\{G\})\,{\cup}\, F_{\HH}(G)$ is $\HH$-nested,
thus any $\GG$-nested set $\Ss$ is an element of $\mathrm{Bl}\,\FF$ as claimed.

Let us now turn to the right-to-left inclusion in (\ref{eq_posets}).\newline
Let $\Ss\,{\in}\,{\rm Bl}_{F_{\HH}(G)}(\FF(\NN(\HH)))$, we have to show 
that~$\Ss$, respectively the set of atoms below~$\Ss$ in $\mathrm{Bl}\,\FF$, 
is nested with respect to $\GG$.

Let us first consider the case when $\Ss$ is $\HH$-nested and does not
contain $F_{\HH}(G)$, i.e., $\Ss$ is one of the elements of the face poset 
$\FF(\NN(\HH))$
that remain after the blowup.  Assume that $\Ss$ is not $\GG$-nested, hence 
there exists an antichain $A$ in $\Ss$ with 
$\bigvee A\,\in\,\GG\,{\setminus}\, \HH$, i.e., $\bigvee A\,{=}\,G$. 
We conclude that $A$ coincides with the set of 
$\HH$-factors of $G$ (cf.~\cite[Prop.\ 2.8.(2)]{FK}),
which contradicts our assumption about $\Ss$ not containing~$F_{\HH}(G)$.

Let us now consider the remaining case, i.e., $\Ss\,{=}\,\Ss'\,{\cup}\,\{G\}$,
where $\Ss'$ is $\HH$-nested, $\Ss'\,{\not \supseteq}\,F_{\HH}(G)$, and
$\Ss'\,{\cup}\,F_{\HH}(G)$ is $\HH$-nested. We have to show that $\Ss$
is $\GG$-nested.

Let $A$ be an antichain  contained  in $\Ss$. If $G\,{\not\in}\, A$, then
$A\,{\subseteq}\,\Ss'$ and $\bigvee A\in \GG\,{\setminus}\, \HH$ implies as
above that $A\,{=}\,F_{\HH}(G)$ contradicting our assumptions. 

If $G\,{\in}\, A$, 
then $A\,{=}\,A'\,{\cup}\,\{G\}$ where $A'$ is an antichain in $\Ss'$.
If $\bigvee A\,{=}\,G$, then $A$ would not be an antichain, hence it suffices to 
show that $\bigvee A\not\in \HH$. 
Consider
\[
    \bigvee A \, \, = \, \,  \bigvee A' \, \, \vee \,G\, \, = \, \, 
        \bigvee A' \,\, \vee \, \,  \bigvee F_{\HH}(G) \, \, = \, \, 
        \bigvee A' \,\,\vee \, \,  
        \bigvee_{{ F\,{\in}\, F_{\HH}(G),\, \,  F \, \mathrm{incom-}    } \atop 
                 {\mathrm{parable\,\, to \,\,elements \, \, 
                                         in\,} A' } 
                } F \, ,                     
\]
where the last equality holds since any $\HH$-factor $F$ of $G$ comparable 
with an element
$a$ in the antichain $A'$ must be smaller than $a$, otherwise $F\,{\geq}\,a$ 
implies $G\,{>}\,a$ which contradicts to $A$ being an antichain. 

We find that $A'\,{\cup}\,\{F{\in} F_{\HH}(G)\,|\, F $ incomparable to elements 
in $A' \}$ is an antichain in $\Ss'\,{\cup}\,F_{\HH}(G)$. With the latter 
being $\HH$-nested by assumption,
we conclude that  $\bigvee A \,{\not \in}\, \HH$ as required, which completes 
our proof. 
\end{pf}

\begin{crl} \label{crl_homeonsc}
Let $\LL$ be a finite atomic meet-semilattice, and $\GG$ a building
set in $\LL$. Then the nested set complex $\NN(\LL,\GG)$ is homeomorphic
to the order complex of~$\LL_{>\hat 0}$,
\[
     \NN(\LL,\GG) \, \, \cong  \, \, \Delta(\LL_{>\hat 0}) \, .
\]
\end{crl}

\begin{pf}
By Theorem~\ref{thm_stsubd} the simplicial fan $\Sigma(\LL_{>\hat 0})$
is a stellar subdivision of $\Sigma(\GG)$ for any building set $\GG$ 
in $\LL$. This in particular implies that the abstract simplicial 
complexes encoding the face structure of the respective fans are 
homeomorphic. The observation that the nested set complex for the 
maximal building set, $\NN(\LL,\LL_{>\hat 0})$, coincides 
with the order complex of $\LL_{>\hat 0}$ finishes our proof. 
\end{pf}

\begin{remrm} \label{rem_top_posets}
  In most of the literature on the topology of posets, the order
  complex of a poset $P$ that has a maximal and a minimal element,
  $\hat 1$ and $\hat 0$, respectively, is understood to be the order
  complex of the proper part 
  $P\,{\setminus}\,\{\hat 0, \hat 1\}$.  Both our theorems can be used
  to study the topology of lattices in this sense:

  Let $\LL$ be a lattice, $\GG$ a building set in $\LL$. We assume that
  $\GG$ contains $\hat 1$, observing that we can always add $\hat 1$
  to a given building set. The nested set complex $\NN(\LL,\GG)$ is a 
  cone with apex $\hat 1$:
\[
  \NN(\LL,\GG)\, \, = \, \, \{\hat 1\} \, \, * \, \, 
                 \NN(\LL,\GG)_{\lceil{\GG\setminus\{\hat 1\}}}\, .
\] 
  Its base can be interpreted as a nested set complex, namely of 
  the meet-semilattice  $\LL\,{\setminus}\,\{\hat1\}$ with respect to 
  $\GG\,{\setminus}\,\{\hat 1\}$:
\[
  \NN(\LL,\GG)_{\lceil{\GG\setminus\{\hat 1\}}} \, \, = \, \, 
      \NN(\LL\,{\setminus}\,\{\hat 1\},\GG\,{\setminus}\,\{\hat 1\})\, .  
\]
 Our theorems state homotopy equivalence, resp.\ homeomorphism, between
 the nested set complex 
 $\NN(\LL\,{\setminus}\,\{\hat 1\},\GG\,{\setminus}\,\{\hat 1\})$ and 
 the order complex of the proper part of $\LL$, 
 $\Delta(\LL\,{\setminus}\,\{\hat 0,\hat 1\})$.  
\end{remrm}




\end{document}